\documentclass[a4paper]{amsart}
\usepackage{amsmath,amsthm,amsfonts,amssymb,amscd,mathrsfs}
\usepackage{psfrag,graphicx}
\usepackage{bbding}
\usepackage{epic,eepic}
\usepackage{color}
\usepackage{ebezier}

\usepackage{epsfig}

\theoremstyle{plain}

\newtheorem{Thm}{Theorem}[section]
\newtheorem{Lem}[Thm]{Lemma}

\theoremstyle{remark}
\newtheorem{Def}[Thm] {Definition}
\newtheorem{Rem}[Thm] {Remark}

\newtheorem{Ex}[Thm] {Example}

\long\def\begcom#1\endcom{}

\newcommand{\length}{\operatorname{\length}}

\def\length{\operatorname{length}}

\begin{document}

\title[Lyapunov-irregularity  ]
{Nonexistence of Lyapunov Exponents for Matrix Cocycles}


\author[X. Tian] {Xueting Tian}
\address[X. Tian]{School of Mathematical Science,  Fudan University\\Shanghai 200433, People's Republic of China}
\email{xuetingtian@fudan.edu.cn}

\keywords{Multiplicative Ergodic Theorem; Lyapunov Exponents; Cocycles;    Specification Property;   Hyperbolic Systems}
\subjclass[2010] { 37H15; 37D20; 37D25; 37C50; 
}
\maketitle

\def\abstractname{\textbf{Abstract}}

\begin{abstract}\addcontentsline{toc}{section}{\bf{English Abstract}}
It follows from Oseledec Multiplicative Ergodic Theorem that the Lyapunov-irregular set of points for which the Oseledec averages of a given continuous cocycle  diverge has zero measure with respect to any   invariant probability measure.    In strong contrast, for   any  dynamical system $f:X\rightarrow X$ with exponential specification property and a H$\ddot{\text{o}}$lder continuous matrix cocycle $A:X\rightarrow G (m,\mathbb{R})$, we show here that if there exist  ergodic measures with different Lyapunov spectrum, then the Lyapunov-irregular set of  $A$ is residual (i.e., containing a dense $G_\delta$ set).


\end{abstract}

\section{Introduction} \setlength{\parindent}{2em}

\subsection{Lyapunov Exponents}



 Let $f$ be an invertible map of a compact metric space $X$ and let $A:X\rightarrow GL (m,\mathbb{R})$  be a  continuous matrix function. One main object of interest is the asymptotic behavior
of the products of $A$ along the orbits of the
transformation $f$, called cocycle induced from $A$: for $n>0$
 $$A (x,n):=A (f^{n-1} (x))\cdots A (f (x))A (x), $$ and $$A (x,-n):=A (f^{-n} (x))^{-1}\cdots A (f^{-2} (x))^{-1} A (f^{-1} (x))^{-1}=A (f^{-n}x,n)^{-1}.$$

 \begin{Def}\label{def2015-Lyapunov-exp}
 For any $x\in X$ and any $0\neq v\in \mathbb{R}^m,$ define the Lyapunov exponent of vector $v$ at $x$,  $$\lambda (A,x,v):=\lim_{n\rightarrow +\infty}\frac1n{\log\|A (x,n)v\|},
   $$ if the limit exists. We say $x$ to be  (forward) {\it Lyapunov-regular} for $A$,  if $\lambda (A,x,v)$ exists for all vector
    $ v\in \mathbb{R}^m\setminus \{0\}.$ Otherwise, $x$ is called to be {\it Lyapunov-irregular} for $A$. Let $LI(A,f)$ denote the space of all Lyapunov-irregular points for $A$.
    \end{Def}

    By Oseledec's Multiplicative Ergodic theorem, for any invariant $\mu$ and $\mu$ a.e. $x$, the Lyapunov exponent
   $\lambda (A,x,v)$ exists at $x$ for all vectors $v\in \mathbb{R}^m.$ In other words, the Lyapunov-irregular set is 
    always   of zero measure for any invariant probability measure. This does not mean that the set of  Lyapunov-irregular points, where the Lyapunov exponents do  not exist, is empty, even if it is completely negligible from the point of view of measure theory. One such interesting result is from \cite{Furman} that    for any uniquely ergodic system, there exists some matrix cocycle  whose Lyapunov-irregular set can be `large' as a set of   second Baire category (also see \cite{Herman,Walter1986,Lenz} for similar discussion).

The notion of  Lyapunov exponent   played important roles in differential dynamical systems, especially in Pesin theory (i.e., Dynamics with non-zero Lyapunov exponents) of {\it Smooth Ergodic Theory}.
It is a quantity that characterizes the rate of separation of infinitesimally close trajectories. The rate of separation can be different for different orientations of initial separation vector. Thus, there is a spectrum of Lyapunov exponents-equal in number to the dimensionality of the phase space.


\bigskip

{\bf Oseledec Multiplicative Ergodic Theorem}  \cite [Theorem 3.4.4]{BP}: Let $f$ be an invertible ergodic measure-preserving transformation of a Lebesgue probability measure space $ (X,\mu).$ Let $A$ be a measurable cocycle whose generator satisfies $\log \|A^\pm (x)\|\in L^1 (X,\mu).$ Then there exist numbers $$\chi_1<\chi_2<\cdots<\chi_l,$$ an $f-$invariant set $\mathcal{R}^\mu$ with $\mu (\mathcal{R}^\mu)=1,$ and an $A-$invariant Lyapunov decomposition of $\mathbb{R}^m$ for $x\in \mathcal{R}^\mu,$
$$\mathbb{R}_x^m=E_{\chi_1} (x)\oplus E_{\chi_2} (x)\oplus \cdots E_{\chi_l} (x)$$ with $dim E_{\chi_i} (x)=m_i,$ such that
for any $i=1,\cdots,l$ and any $0\neq v\in E_{\chi_i} (x)$ one has
$$\lim_{n\rightarrow \pm\infty} \frac1n \log \|A (x,n)v\|=\chi_i$$
and $$\lim_{n\rightarrow \pm\infty} \frac1n \log det A (x,n) =\sum_{i=1}^lm_i\chi_i.$$

 \begin{Def}\label{Def-LyapunovExponents}
The numbers $\chi_1,\chi_2,\cdots,\chi_l$ are called the {\it Lyapunov exponents} of measure $\mu$ for cocycle $A$ and the dimension $m_i$ of the space $E_{\chi_i} (x)$ is called the {\it multiplicity} of the exponent $\chi_i.$ The collection of pairs $$Sp(\mu,A)=\{(\chi_i,m_i):1\leq i \leq l\}$$ is the {\it Lyapunov spectrum} of measure $\mu.$
\end{Def}

Remark that for any ergodic measure $\mu,$ all the points in the set $\mathcal{R}^\mu$ are  Lyapunov-regular.

\subsection{Results}

Recall that $Y$ is called residual in $X,$ if $Y$ contains a dense $G_\delta$ subset of $X.$ The notion of residual set is usually used to describe a set being `large' in the
topological sense.

 \begin{Thm}\label{Thm-RESIDUAL-LyapunovIrregular-HolderCocycle} 
Let $f:X\rightarrow X$ be a     homeomorphism of a compact metric space $X$ with exponential  specification.     Let $A:X\rightarrow GL (m,\mathbb{R})$  be a H$\ddot{o}$der  continuous matrix function. Then either all ergodic measures have same Lyapunov spectrum or the  Lyapunov-irregular set $LI (A,f)$ is  residual in $X.$
\end{Thm}

\begin{Rem}  From \cite [Theorem 4]{Furman} we know    for any uniquely ergodic system, there exists some continuous matrix cocycle  whose Lyapunov-irregular set is a dense set of   second Baire category. So, it is admissible  to satisfy  that all ergodic measures have same Lyapunov spectrum and  simultaneously   Lyapunov-irregular points form a set of  second Baire category. On the other hand, if $m=1,$
the above phenomena of \cite{Furman} naturally does not happen for any dynamical system. Let us explain more precisely. If $m=1,$ the Lyapunov exponent can be written as Birkhoff ergodic average $$ \lim_{n\rightarrow +\infty}\frac1n\sum_{j=0}^{n-1}\phi(f^j(x))$$ where $\phi(x)=\log\|A (x)\|$ is a continuous function.  If all ergodic measures have same Lyapunov spectrum, then by Ergodic Decomposition theorem so do  all invariant measures and thus by weak$^*$ topology, the limit  $  \lim_{n\rightarrow +\infty}\frac1n\sum_{j=0}^{n-1}\phi(f^j(x))$  should exist at every point $x\in X$ and equal to the given spectrum. Moreover,   the case of $m=1$ is in fact to study Birkhoff ergodic average and it has been studied for systems with specification or its variants by many authors,   see \cite{BarLV,LW2,LW,APT,Olsen} and reference therein.
\end{Rem}

As a particular case of Theorem \ref{Thm-RESIDUAL-LyapunovIrregular-HolderCocycle} we have a consequence for the derivative cocycle of      hyperbolic systems.
 Let $LI(f):=LI(Df,f)$. It is called Lyapunov-irregular set of system $f$.

\begin{Thm}\label{Thm-1Horseshoe-RESIDUAL-LyapunovIrregular-Cocycle}
 Let $f$ be  a $C^{1+\alpha}$   diffeomorphism of a compact Riemanian manifold $M$
 and $X\subseteq M$ be a topologically mixing locally maximal invariant subset. Then either all ergodic measures have same Lyapunov spectrum or the
Lyapunov-irregular set $LI (f)$ 
is residual in $X.$
\end{Thm}

As said in \cite{Bar-Gel} 
that the study of Lyapunov exponents lacks today a satisfactory general approach for non-conformal maps, since a complete understanding is just known for some cases such as requiring a clear separation of Lyapunov directions or some number-theoretical properties etc, see \cite{Bar-Gel} and its references therein.   In  this program our Theorem \ref{Thm-1Horseshoe-RESIDUAL-LyapunovIrregular-Cocycle} gives one such characterization on Lyapunov exponents for general non-conformal hyperbolic dynamics.




\section{Specification and Lyapunov Metric}

\subsection{Specification and Exponential Specification}


Now we introduce  {\it (exponential) specification} property. Let $f$ be a continuous map of a compact metric space $X$.

\begin{Def}\label{Def-Specification}
$f$ is called to have  specification property, if the
following holds: for any $\delta>0$ there exists an increasing
sequence of integers $N=N (\delta)>0$ such that for any
$k\geq 1 $, any $k$ points $x_1,x_2,\cdots,x_k\in X$,
 any integers $a_1\leq b_1<a_2\leq b_2<\cdots<a_k\leq b_k$ with $a_{i+1}-b_i\geq N,$  there exists a point $y\in
 X$
such that  $d (f^{i} (y),
f^i (x_j))<\delta,\,\,a_j\leq i\leq b_j, \,\, 1\leq j\leq k.$

\end{Def}

Remark that the specification property introduced by Bowen\cite{Bowen3} required that the shadowing point $y$ is periodic. That is, for any $p\geq b_k-a_1+N,$ the chosen point $y$ in above definition further
satisfies $f^p (y)=y.$ We call this to be   {\it Bowen's   Specification property.}

\begin{Def}\label{Def-exp-Specification}
  $f$ is called to have  exponential specification property with exponent $\lambda>0$  (only dependent on the system $f$ itself), if specification property holds and the   inequality in specification can be shadowed exponentially, i.e., $$d (f^i (x_j),f^i (y))<\delta e^{-\lambda\min\{i-a_j,b_j-i\}},\,\,a_j\leq i\leq b_j, \,\, 1\leq j\leq k.$$

  If further the tracing point $y$ is periodic with $f^p(y)=y$, then we say $f$ has Bowen's exponential specification property.

\end{Def}

For convenience, we say the orbit segments $x,fx,\cdots,f^nx$ and $y,fy,\cdots,f^ny$ are exponentially $\delta$ close with exponent $\lambda$, meaning that $$d (f^i (x),f^i (y))<\delta e^{-\lambda\min\{i,n-i\}},0\leq i\leq n.$$

\begin{Rem}\label{Rem-exp-spe}
It is not difficult to see that
  $$\text{   (Bowen's) Specification
+ Local product structure}$$ $$
\Rightarrow \text{  (Bowen's)  Exponential  specification.}$$
Recall that every hyperbolic set has local product structure and every topologically mixing locally maximal hyperbolic set has Bowen's specification property\cite{Bowen3}.
So
every topologically mixing locally maximal hyperbolic set has exponential specification property.
   As a particular case,  every transitive Anosov diffeomorphism  has    exponential
specification property, since it is known that every transitive Anosov diffeomorphism   is  topologically mixing. If  a homeomorphism $f$ is topologically conjugated to a homeomorphism $g$ satisfying (Bowen's)  exponential specification property  with exponent for
 some $\beta>0$, and the inverse conjugation is $\gamma-$H$\ddot{\text{o}}$lder continuous, then it is not difficult to see that $f$   has  (Bowen's) exponential specification property  with exponent $\beta\gamma>0.$
\end{Rem}

\begin{Ex}\label{exm-nonhyperbolic} {\it Some non-hyperbolic systems with exponential specification:

(1) From \cite{Go} we know that  non-hyperbolic diffemorphism $f$ with $C^{1+Lip}$ smoothness, conjugated to a transitive Anosov diffeomorphism, exists even the conjugation and its inverse  is H$\ddot{\text{o}}$lder continuous. This example satisfies Bowen's  exponential specification property.

  (2)  For time-1 map  of a geodesic
flow of compact connected negative curvature manifolds, it is a partially  (non-hyperbolic) hyperbolic dynamical system. 
 Its  exponential  specification property   can be deduced from the local product structure and specification property of the flow which is naturally hyperbolic, see \cite{Bow}.  Here we remark that  the shadowing point   may be not periodic with respect to the  time-1 map  because the shadowing of flow has a small time reparameterization.  }


\end{Ex}





\subsection{Lyapunov Exponents and Lyapunov Metric}

Suppose  $f:X\rightarrow X$ to be an invertible map on a compact metric space $X$ and $A:X\rightarrow GL (m,\mathbb{R})$ to  be a  continuous matrix function.
We denote the standard scalar product in $\mathbb{R}^m$ by $<\cdot,\cdot>$. For a fixed $\epsilon>0$ and a regular point $x$ we introduce the {\it $\epsilon-$Lyapunov scalar product  (or metric)} $<\cdot,\cdot>_{x,\epsilon}$ in $\mathbb{R}^m$ as follows. For $u\in E_{\chi_i} (x),\, v\in E_{\chi_j} (x),\,i\neq j$ we define $<\cdot,\cdot>_{x,\epsilon}=0.$ For $i=1,\cdots,l$ and $u,v\in u\in E_{\chi_i} (x),$ we define
$$<\cdot,\cdot>_{x,\epsilon}=m\sum_{n\in\mathbb{Z}}<A (x,n)u,A (x,n)v>exp (-2\chi_in-\epsilon |n|).$$
Note that the series converges exponentially for any regular $x$. The constant $m$ in front of the conventional formula is introduced for more convenient comparison with the standard scalar product. Usually, $\epsilon$ will be fixed and we will denote $<\cdot,\cdot>_{x,\epsilon}$ simply by $<\cdot,\cdot>_{x}$
 and call it the {\it Lyapunov scalar product.} The norm generated by this scalar product is called the {\it Lyapunov norm}  and is denoted by $\|\cdot\|_{x,\epsilon}$ or $\|\cdot\|_{x}.$

 We summarize below some important properties of the Lyapunov scalar product and norm; for more details  see  \cite [\S 3.5.1-3.5.3]{BP}.  A direct calculation shows   \cite [Theorem 3.5.5]{BP} that for any regular $x$ and any $u\in E_{\chi_i} (x)$
 \begin{eqnarray}\label{Lyapunov-norm} \,\,\,\,\,\,\,\,
  exp (n\chi_i-\epsilon|n|)\|u\|_{x,\epsilon}\leq \|A (x,n)u\|_{f^nx,\epsilon}\leq exp  (n\chi_i+\epsilon|n|)\|u\|_{x,\epsilon}\,\,\,\,\,\forall n\in \mathbb{Z},
 \end{eqnarray}
\begin{eqnarray}\label{Lyapunov-norm-Maximal-ightarrow}
 exp (n\chi-\epsilon|n|) \leq \|A (x,n)u\|_{f^nx\leftarrow x}\leq exp  (n\chi+\epsilon|n|) \,\,\,\,\,\forall n\in \mathbb{Z},
 \end{eqnarray}
where $\chi=\chi_l$  is the maximal Lyapunov exponent and $\|\cdot\|_{f^nx\leftarrow x}$ is the operator norm with respect to the Lyapunov norms. It is defined for any matrix $A$ and any regular points $x,y$ as follows:
$$\|A\|_{y\leftarrow x}=\sup\{\|Au\|_{y,\epsilon}\cdot \|u\|^{-1}_{x,\epsilon}:0\neq u\in \mathbb{R}^m\}.$$

We emphasize that, for any given $\epsilon>0,$ Lyapunov scalar product and Lyapunov norm are defined only for regular points with respect to the given measure. They depend only measurably on the point even if the cocycle is H$\ddot{\text{o}}$lder. Therefore, comparison with the standard norm becomes important. The uniform lower bound follow easily from the definition: $$\|u\|_{x,\epsilon}\geq \|u\|.$$ The upper bound is not uniform, but it changes slowly along the regular orbits (\cite{BP}, Prop. 3.5.8): there exists a measurable function $K_\epsilon (x)$ defined on the set of regular points $\mathcal{R}^\mu$ such that
\begin{eqnarray}\label{eq-different-norm-estimate}
\|u\|\leq \|u\|_{x,\epsilon}\leq K_\epsilon (x)\|u\|\,\,\,\,\,\,\,\,\,\forall x\in\mathcal{R}^\mu,\,\,\forall u\in  \mathbb{R}^m
 \end{eqnarray}
\begin{eqnarray}\label{eq-estimate-K-epsilon}
  K_\epsilon (x)e^{-\epsilon n}\leq K_\epsilon (f^nx)  \leq K_\epsilon (x)e^{\epsilon n}\,\,\,\,\,\,\,\,\,\forall x\in\mathcal{R}^\mu,\,\,\forall n\in  \mathbb{Z}.
 \end{eqnarray}
 These estimates are obtained in \cite{BP} using the fact that $\|u\|_{x,\epsilon}$ is {\it tempered}, but they can also be checked directly using the definition of  $\|u\|_{x,\epsilon}$ on each Lyapunov space and noting that angles between the spaces change slowly.

 For any matrix $A$ and any regular points $x,y,$ inequalities  (\ref{eq-different-norm-estimate}) and  (\ref{eq-estimate-K-epsilon}) yield
\begin{eqnarray}\label{eq-estimate-norm-K-epsilon}
  K_\epsilon (x)^{-1}\|A\| \leq \|A\|_{y\leftarrow x}\leq K_\epsilon (y) \|A\|.
 \end{eqnarray}

When $\epsilon$ is fixed we will usually omit it and write $K (x)=K_\epsilon (x).$ For any $l>1$ we also define the following sets of regular points
\begin{eqnarray}\label{eq-estimate-measure-Pesinblock}
  \mathcal{R}^\mu_{\epsilon,l}=\{x\in \mathcal{R}^\mu: \,\,K_\epsilon (x)\leq l\}.
 \end{eqnarray}
Note that $\mu (\mathcal{R}^\mu_{\epsilon,l} )\rightarrow 1$ as $l\rightarrow \infty.$ Without loss of generality, we can assume that the set $\mathcal{R}^\mu_{\epsilon,l}$ is compact and that Lyapunov splitting and Lyapunov scalar product are continuous on $\mathcal{R}^\mu_{\epsilon,l}.$ Indeed, by Luzin's theorem we can always find a subset of $\mathcal{R}^\mu_{\epsilon,l}$ satisfying these properties with arbitrarily small loss of measure (for standard Pesin sets these properties are automatically satisfied).

\section{Norm  Estimate of Cocycles \& Generic Property of  Lyapunov-irregularity}


\subsection{Estimate of the norm of H$\ddot{\text{o}}$lder cocycles}\label{Estimate-norm}

Before proving  generic property of  Lyapunov-irregularity, we need to recall some useful lemmas   as follows. Firstly let us recall a general estimate of the norm of $A$ along any orbit segment close to a regular one\cite{Kal}.

 \begin{Lem}\label{Lem-simple-estimate-Lyapunov}   \cite [Lemma 3.1] {Kal}
Let $A$ be an $\alpha-$H$\ddot{\text{o}}$lder cocycle ($\alpha>0$) over a continuous map $f$ of a compact metric space $X$ and let $\mu$ be an ergodic measure for $f$ with the largest Lypunov exponent $\chi.$ Then for any positive $\lambda$ and $\epsilon$ satisfying $\lambda>\epsilon/ \alpha$ there exists $c>0$ such that for any $n\in\mathbb{N}$, any regular point $x$ with both $x$ and $f^nx$ in $\mathcal{R}^\mu_{\epsilon,l}$, and any point $y\in X$ such that the orbit segments $x,fx,\cdots,f^nx$ and $y,fy,\cdots,f^n (y)$ are exponentially $\delta$ close with exponent $\lambda$ for some $\delta>0$ we have
\begin{eqnarray}\label{eq-estimate-simple-Lyapunov-1}
  \|A (y,n)\|_{f^nx\leftarrow x}\leq e^{cl\delta^\alpha}e^{n (\chi+\epsilon)}\leq e^{2n\epsilon+cl\delta^\alpha}\|A (x,n)\|_{f^nx\leftarrow x}
 \end{eqnarray}
 and
 \begin{eqnarray}\label{eq-estimate-simple-Lyapunov-2}
  \|A (y,n)\| \leq l^2 e^{cl\delta^\alpha}e^{n (\chi+\epsilon)}\leq l^2 e^{2n\epsilon+cl\delta^\alpha}\|A (x,n)\|.
 \end{eqnarray}
 The constant $c$ depends only on the cocycle $A$ and on the number $ (\alpha\lambda-\epsilon).$

\end{Lem}

\begin{Lem}\label{Lem-New-simple-estimate-Lyapunov-new}
Let $A$ be an $\alpha-$H$\ddot{\text{o}}$lder cocycle  ($\alpha>0$) over a continuous map $f$ of a compact metric space $X$ and let $\mu$ be an ergodic measure for $f$ with the largest Lypunov exponent $\chi.$ Then for any positive $\lambda$ and $\epsilon$ satisfying $\lambda>\epsilon/ \alpha$ there exists $\delta>0$ such that for any $n\in\mathbb{N}$, any regular point $x$ with both $x$ and $f^nx$ in $\mathcal{R}^\mu_{\epsilon,l}$, and any point $y\in X$ such that the orbit segments $x,fx,\cdots,f^nx$ and $y,fy,\cdots,f^n (y)$ are exponentially $\delta$ close with exponent $\lambda$ for some $\delta>0$ we have
 \begin{eqnarray}\label{eq-estimate-simple-Lyapunov-new}
  \|A (y,n)\| \leq l^2 e^{l}e^{n (\chi+\epsilon)}\leq l^2e^l e^{2n\epsilon} \|A (x,n)\|.
 \end{eqnarray}

\end{Lem}

{\bf Proof.} For Lemma \ref{Lem-simple-estimate-Lyapunov}, let $\delta>0$ small enough such that $$c\delta^\alpha<1.$$ Then the   estimate (\ref{eq-estimate-simple-Lyapunov-new}) is obvious from Lemma \ref{Lem-simple-estimate-Lyapunov}. \qed

\bigskip

Another lemma is to estimate the growth of vectors in a ceratin cone $K\subseteq \mathbb{R}^m$ invariant under $A (x,n)$ \cite{Kal}. Let $x$ be a point in $\mathcal{R}^\mu_{\epsilon,l} $ and $y\in X$ be a point such that the orbit segments $x,fx,\cdots,f^nx$ and $y,fy,\cdots,f^ny$ are exponentially $\delta$ close with exponent $\lambda.$ We denote
$x_i=f^ix$ and $y_i=f^iy,\,i=0,1,\cdots,n.$ For each $i$ we have orthogonal splitting $\mathbb{R}^m=E_i\oplus F_i,$ where $E_i$ is the Lyapunov space at $x_i$ corresponding to the largest Lyapunov exponent $\chi$ and $F_i$ is the direct sum of all other Lyapunov spaces at $x_i$ corresponding to the Lyapunov exponents less than $\chi.$ For any vector $u\in \mathbb{R}^m$ we denote by $u=u'+u^\perp$ the corresponding splitting with $u'\in E_i$ and $u^\perp\in F_i;$ the choice of $i$ will be clear from the context. To simplify notation, we write $\|\cdot\|_i$ for the Lyapunov norm at $x_i$. For each $i=0,1,\cdots,n$ we consider cones $$K_i=\{u\in \mathbb{R}^m:\,\|u^\bot\|_i\leq \|u'\|_i\}\,\,\,\,\text{ and } \,\,\,\,K_i^\eta=\{u\in \mathbb{R}^m:\,\|u^\bot \|_i\leq  (1-\eta)\|u'\|_i\}$$
with $\eta>0$. Remark that \begin{eqnarray}\label{eq-u-u'}
\|u\|_i\geq \|u'\|_i\geq\frac1{\sqrt2}\|u\|_i.
\end{eqnarray}
We will consider the case when $\chi$ is not the only Lyapunov exponent of $A$ with respcet to $\mu.$ Otherwise $F_i=\{0\}, K_i^\eta=K_i=\mathbb{R}^m$, and the argument becomes simpler. Recall that $\epsilon<\epsilon_0=\min\{\lambda\alpha, (\chi-\nu)/2\},$ where $\nu<\chi$ is the second largest Lyapunov exponent of $A$ with respect to $\mu.$

\begin{Lem}\label{Lem-simple-estimate-Lyapunov-2}  \cite [Lemma 3.3]{Kal}
In the notation above, for any regular set $\mathcal{R}^\mu_{\epsilon,l}$, there exist $\eta,\delta>0$ such that if  $x,f^nx\in\mathcal{R}^\mu_{\epsilon,l}$  and   the orbit segments $x,fx,\cdots,f^nx$ and $y,fy,\cdots,f^n (y)$ are exponentially $\delta$ close with exponent $\lambda$, then for every $i=0,1,\cdots,n-1$ we have $A (y_i) (K_i)\subseteq K_i^\eta$ and  $\| (A (y_i)u)'\|_{i+1}\geq e^{\chi-2\epsilon}\|u'\|_i$ for any $u\in K_i.$

\end{Lem}

\subsection{Residual Property of Maximal Lyapunov-irregularity}


The Maximal Lyapunov exponent  (or simply, MLE) of $A:X\rightarrow GL (m,\mathbb{R})$ at one point $x\in X$ is defined as the limit $$\chi_{max} (A,x):=\lim_{n\rightarrow \infty}\frac1n{\log\|A (x,n)\|},
   $$ if it exists. In this case $x$ is called to be (forward) {\it Max-Lyapunov-regular}. 
    Otherwise, $x$ is {\it Max-Lyapunov-irregular}. By Oseledec's Multiplicative Ergodic Theorem  (or Kingman's Sub-additional Ergodic Theorem), for
    any ergodic  measure $\mu$ and $\mu$ a.e. point $x$,  MLE always exists and is constant, denoted by $\chi_{max}(A,\mu)$. 
Let   $MLI (A,f)$ denote the set of all   Max-Lyapunov-irregular points.
 Then it is of zero measure for any ergodic measure and  by Ergodic Decomposition theorem so does it for all invariant measures.
 Now let us show a residual result for $MLI (A,f)$.




 \begin{Thm}\label{Thm-Maximal-Residual-LyapunovIrregular-Cocycle}

Let $f:X\rightarrow X$ be a     continuous map of a compact metric space $X$ with exponential  specification.     Let $A:X\rightarrow GL (m,\mathbb{R})$  be a H$\ddot{o}$der continuous function. Suppose that $$\inf_{\mu\in \mathcal{M}_{erg} (X)}\chi_{max} (A,\mu)<\sup_{\mu\in \mathcal{M}_{erg} (X)}\chi_{max} (A,\mu).$$ Then    the  Max-Lyapunov-irregular set $MLI (A,f)$  is   residual in $X$.
\end{Thm}

In other words, either all ergodic measures have same maximal Lyapunov exponent  or Max-Lyapunov-irregular set $MLI (A,f)$  is   residual in $X$.

{\bf Proof. }    Take two ergodic measures $\mu$ and $\omega$ such that $$\int \chi_{max} (A,x)d\mu>\int \chi_{max} (A,x)d\omega.$$ If let $a=\int \chi_{max} (A,x)d\mu$ and $b=\int \chi_{max} (A,x)d\omega$, we can choose  $\tau>0$ such that  $a-2\tau>b+2\tau.$ Let $C=\max_{x\in X}\{\|A(x)\|,\|A^{-1}(x)\|\}.$

Define $$O_n:=\{w |\,\,\exists\,\,n_1,n_2>n\,\,s.t.\,\,\,\frac1{n_1}\log \|A (w,{n_1})\|>a-\tau\,\,\&\,\,\frac1{n_2}\log \|A (w,{n_2})\|<b+\tau\}.$$ By continuity of $A (x,n)$, $O_n$ is open. It is straightforward to check that $$\bigcap_{n\geq 1}O_n\subseteq MLI (A,f).$$ So we only need to prove that $O_n$ is dense in $X.$ Fix $x_0\in X$ and $t>0$, we will show $O_n\cap B(x_0,t)\neq \emptyset.$

More precisely, fix $n\geq 1$ and  firstly let us recall or define some notations. Let $A$ be   $\alpha-$H$\ddot{\text{o}}$lder continuous and let $C:=\max_{x\in X}\|A^\pm (x)\|$. Let $\lambda $ be the positive number in the definition of exponential specification. Take  $\epsilon\in (0,\frac12\tau)$  satisfying $\lambda>\epsilon/ \alpha$.
 For the measures $\mu$ and $\omega,$ take $l$ large enough such that $$\mu (\mathcal{R}^\mu_{\epsilon,l})>0,\,\,\,\omega (\mathcal{R}^\omega_{\epsilon,l})>0.$$

 Take $\eta>0,\delta\in(0,t)$ small enough such that it is applicable to Lemma \ref{Lem-New-simple-estimate-Lyapunov-new} and Lemma \ref{Lem-simple-estimate-Lyapunov-2}. For $\delta,\,\lambda,$ there is $N=N$ satisfies the exponential specification.

By Poincar$\acute{\text{e}}$ Recurrence theorem, there exist two points $x\in \mathcal{R}^\mu_{\epsilon,l},\,\,z\in \mathcal{R}^\omega_{\epsilon,l}$ and two increasing sequences $\{H_i\},\,\,\{L_i\}\nearrow \infty$ such that
$f^{H_i} (x)\in \mathcal{R}^\mu_{\epsilon,l},\,\,f^{L_i} (z)\in \mathcal{R}^\omega_{\epsilon,l}.$
Take $H=H_i\gg \max\{ N,n\}$ such that $$\frac1{\sqrt 2\,l} e^{H (a-2\epsilon)}>C^N e ^{(H+N) (a-\tau)}$$ and take $L=L_j\gg H+N$ large enough such that $$ l^2 e^{l}e^{L (b+\epsilon)}C^{H+2N}<e^{ (b+\tau) (L+H+2N)}.$$

Now let us consider three orbit segments $$\{x_0\},\,\,\{x,fx,\cdots,f^{H}x\},\,\,\{z,fz,\cdots,f^{L}z\}$$ (hint: $a_1=b_1=0,a_2=N,b_2=a_2+H,a_3=b_2+N,b_3=a_3+L$)  for the exponential specification. Then there is $y_0\in X$ such that $d(y_0,x_0)<\delta ,$ the orbit segments $x,fx,\cdots,f^{H}x$ and $y,fy,\cdots,f^Hy$ are exponentially $\delta$ close with exponent $\lambda$ where $y=f^{N}y_0$, and simultaneously the orbit segments $z,fz,\cdots,f^{L}z$ and $y',fy',\cdots,f^Ly'$ are exponentially $\delta$ close with exponent $\lambda$ where $y'=f^{H+2N}y_0$.

Firstly let us consider the orbit segments $x,fx,\cdots,f^{H}x$ and $y,fy,\cdots,f^Hy$.
By   Lemma \ref{Lem-simple-estimate-Lyapunov-2} (in this estimate $\chi=a$, being the largest Lyapunov exponent of $\omega$),   for any $u\in K_0$ with $\|u\|=1,$
$$\| (A (y,H)u)'\|_H \geq  e^{H (a-2\epsilon)}\|u'\|_0.$$ 
 Together with (\ref{eq-different-norm-estimate}) and (\ref{eq-u-u'}),  we have
 $$\|A (y,H)\|\geq\|A (y,H)u\| \geq\frac1l\|A (y,H)u\|_H \geq\frac1l\| (A (y,H)u)'\|_H \geq \frac1l e^{H (a-2\epsilon)}\|u'\|_0$$
$$\geq  \frac 1{\sqrt 2\,l}e^{H (a-2\epsilon)}\|u\|_0\geq \frac 1{\sqrt 2\,l}e^{H (a-2\epsilon)}\|u\|=\frac1{\sqrt 2\,l} e^{H (a-2\epsilon)}>C^{N}e ^{(H+N) (a-\tau)} .$$
Then \begin{eqnarray}\label{eq-2015-1} & &\|A (y_0,H+N)\|\geq m(A(y_0,N)) \cdot \|A (y,H)\| \\
&\geq& C^{-N} \cdot \|A (y,H)\|>e ^{(H +N) (a-\tau)},\nonumber
\end{eqnarray}
 where
 $m(B)$ denotes the minimal norm of linear map $B$.

Secondly let us consider the orbit segments $z,fz,\cdots,f^{L}z$ and $y',fy',\cdots,f^Ly'$.
By the first estimate in (\ref{eq-estimate-simple-Lyapunov-new})  of Lemma \ref{Lem-New-simple-estimate-Lyapunov-new} (in this estimate $\chi=b$, being  the largest Lyapunov exponent of $\mu$)   we have
$$\|A (y',L)\|
\leq  l^2 e^{l}e^{L (b+\epsilon)}. $$
Then
\begin{eqnarray}\label{eq-2015-2}& &\|A (y_0,L+H+2N)\|\leq \|A (y',L)\|\cdot \|A (y_0,H+2N)\| \\
 &\leq&  l^2 e^{l}e^{L (b+\epsilon)}C^{H+2N}<e^{ (b+\tau) (L+H+2N)}. \nonumber
\end{eqnarray}

Take $n_1=H+N$ and $n_2=L+H+2N$, then (\ref{eq-2015-1}) and (\ref{eq-2015-2}) imply $y_0\in O_n .$  Recall $d(y_0,x_0)<\delta$ and $ \delta<t$ so that $y_0\in   B(x_0,t).$ So we complete the proof. \qed

\subsection{Proof of Theorem \ref{Thm-RESIDUAL-LyapunovIrregular-HolderCocycle} and \ref{Thm-1Horseshoe-RESIDUAL-LyapunovIrregular-Cocycle}}

For    a Lyapunov-regular point $x$ of cocycle  $A$, let $\lambda_1(x)\geq \lambda_2(x)\geq\cdots \geq \lambda_m(x)$ denote the Lyapunov exponents of $x$ for $A$.
 Let $$\Lambda^A_i(x)=\sum_{j=1}^i\lambda_j(x).
 $$ and for any invariant measure $\mu,$ define
 $ \Lambda^A_i(\mu)=\int \Lambda^A_i(x)d\mu.$ 
Then it is easy to compute  to obtain that: for any two ergodic measures $\mu,\nu\in \mathcal{M}^e_{f}(X)$,
\begin{eqnarray}\label{Equivalent-Lyapunov-exponents}
Sp(\mu,A)=Sp(\nu,A)\Leftrightarrow \Lambda^A_i(\mu)=\Lambda^A_i(\nu),\,\forall \,i.
\end{eqnarray}
Let us consider cocycle $\wedge^i A(x,n)$ induced by cocycle $A(x,n)$ on the $i$-fold exterior powers $\wedge^i \mathbb{R}^m$.
 For any $1\leq i\leq m$ and Lyapunov-regular point $x \in X,$
\begin{eqnarray}\label{Equivalent-Lyapunov-exponents-2015} \lim_{n\rightarrow \infty}\frac1n\log \|\wedge^i A(x,n)\|=\sum_{j=1}^i\lambda_j(x)=\Lambda^A_i(x).
\end{eqnarray}
 This trick is related to Ragunatan's proof of the Multiplicative
Ergodic Theorem \cite [\S 3.4.4] {BP} and was also used in \cite{Kal,WS}.
\bigskip

{\bf Proof of Theorem \ref{Thm-RESIDUAL-LyapunovIrregular-HolderCocycle}.}
Assume that there are two ergodic measures with different Lyapunov spectrum. By  (\ref{Equivalent-Lyapunov-exponents}) and (\ref{Equivalent-Lyapunov-exponents-2015}), there is some $1\leq i \leq m$ such that
\begin{eqnarray}\label{eq-twodifferent-measure}
\inf_{\mu\in \mathcal{M}^e_{f}(X)}\lim_{n\rightarrow \infty}\frac1n\log
\|\wedge^i A(x,n)\|< \sup_{\mu\in \mathcal{M}^e_{f}(X)}\lim_{n\rightarrow \infty}\frac1n\log \|\wedge^i A(x,n)\|.
\end{eqnarray}
 Then we can apply Theorem \ref{Thm-Maximal-Residual-LyapunovIrregular-Cocycle} to the cocycle
$\wedge^i A(x,n)$ and obtain that the Max-Lyapunov-irregular set of $\wedge^i A$ $ MLI(\wedge^i A,f)$ is residual in $X$. Note that $LI(A,f)\supseteq   MLI(\wedge^i A,f)$, since by (\ref{Equivalent-Lyapunov-exponents-2015}) a point Lyapunov-regular for $A$ is also Lypunov-regular for $\wedge^i A$.  So $LI(A,f)$  is also residual in $X$. Now we complete the proof.
\qed
\bigskip

{\bf Proof of Theorem \ref{Thm-1Horseshoe-RESIDUAL-LyapunovIrregular-Cocycle}} From Remark \ref{Rem-exp-spe} we know  $f|_X$ has exponential specification. Applying Theorem \ref{Thm-RESIDUAL-LyapunovIrregular-HolderCocycle} for cocycle $A(x,n)=D_xf^n$, one ends the proof. \qed

\section*{Acknowlegements}   The research of X. Tian was  supported by National Natural Science Foundation of China  (grant no. 11301088) and  Specialized
  Research Fund for the Doctoral Program of Higher Education  (No.  20130071120026).



\begin{thebibliography}{10}


\itemsep -2pt

\small
\bibitem{APT} S. Albeverio, M. Pratsiovytyi and G. Torbin, {\it Topological and fractal properties of subsets of
real numbers which are not normal,}  Bull. Sci. Math., 129  (2005), 615-630.

\bibitem{Olsen2}I.-S. Baek and L. Olsen, {\it Baire category and extremely non-normal points of invariant sets of
IFS's,}  Discrete Contin. Dyn. Syst., 27  (2010), 935-943.






\bibitem{Bar-Gel}
 L. Barreira, K. Gelfert, {\it Dimension estimates in smooth dynamics: a survey of recent results},  Ergodic Theory and Dynamical Systems, 2011, 31(03): 641-671.


\bibitem{BarLV} L. Barreira, J. Li , C. Valls,  {\it Irregular sets for ratios of Birkhoff averages are residual},  Publicacions Matem$\grave{\text{a}}$tiques, 2014, 58: 49-62.



\bibitem{BP} Luis Barreira and Yakov B. Pesin, {\it  Nonuniform hyperbolicity}, Cambridge Univ. Press, Cambridge  (2007).








\bibitem{Bow}R. Bowen, Periodic orbits for hyperbolic flows, Amer. J. Math., 94 (1972), 1-30.

\bibitem{Bowen3} R. Bowen, Periodic points and measures for Axiom A diffeomorphisms, Trans. Amer. Math. Soc. 154 (1971), 377-397.


















\bibitem{Furman}  A. Furman, {\it  On the multiplicative ergodic theorem for uniquely ergodic systems}, Annales de l'Institut Henri Poincar$\acute{\text{e}}$ (B) Probability and Statistics, 1997, 33(6): 797-815.


\bibitem{Hyde} J. Hyde, V. Laschos, L. Olsen, I. Petrykiewicz and A. Shaw, {\it Iterated Cesaro averages, fre-
quencies of digits and Baire category,} Acta Arith., 144  (2010), 287-293.



\bibitem{Go} A. Gogolev, {\it Diffeomorphisms H$\ddot{\text{o}}$lder conjugate to Anosov diffeomorphisms,}  Ergodic Theory and Dynamical Systems, Vol. 30, no. 2  (2010),  441-456




\bibitem{Herman} M. Herman, {\it Construction d'un diff\'{e}omorphisme minimal d'entropie topologique non nulle,} Ergod.
Th. Dynam. Sys. 1 (1981), 65-76.






\bibitem{Kal} B. Kalinin,  {\it Liv$\check{s}$ic Theorem for matrix cocycles,} Annals of mathematics, 2011, 173 (2), 1025-1042.



\bibitem{Lenz} D.  Lenz, {\it Existence of non-uniform cocycles on uniquely ergodic systems}, Annales de l'Institut Henri Poincar$\acute{\text{e}}$ (B) Probability and Statistics, 2004, 40(2): 197-206.

\bibitem{LW2}  J.   Li and M. Wu, {\it The sets of divergence points of self-similar measures are residual,}  J.
Math. Anal. Appl., 404  (2013), 429-437.

\bibitem{LW}   J. Li, M. Wu, {\it Generic property of irregular sets in systems satisfying the specificaiton property,} Discrete and Continuous Dynamical Systems 34 (2014), 635-645.







\bibitem{Olsen} L. Olsen, {\it Extremely non-normal numbers,}  Math. Proc. Cambridge Philos. Soc., 137  (2004),
43-53.


\bibitem{Os} V. I. Oseledec, {\it Multiplicative ergodic theorem, Liapunov
characteristic numbers for dynamical systems}, Trans. Moscow Math.
Soc., 19  (1968), 197-221; translated from Russian.


























\bibitem{Walter1986}  P. Walters, {\it Unique ergodicity and random matrix products,}  Lyapunov Exponents, Springer Berlin Heidelberg, 1986: 37-55.

\bibitem{Walter} P. Walters,  {\it An introduction to ergodic theory,}
Springer-Verlag, 2001.

\bibitem{WS}Z. Wang and W. Sun, {\it Lyapunov exponents of hyperbolic measures and hyperbolic period
orbits,} Trans. Amer. math. Soc., 362   (2010), 4267-4282.












\end{thebibliography}
\end{document}